\documentclass[11pt,twoside,drfat]{article}
\usepackage{amssymb,amsthm}
\usepackage[namelimits,sumlimits,fleqn]{amsmath}
\usepackage{setspace}
\usepackage[ hmargin={2.2cm,2.4cm}, top=36mm,
    a4paper,pdftex, textheight=230mm,
    headheight=20mm, headsep=10mm, bindingoffset=0.9cm]
    {geometry}
\usepackage[small, margin=20pt]{caption}
\usepackage[usenames,dvipsnames]{color}
\definecolor{darkblue}{rgb}{0.0,0.0,0.6}
\usepackage{float}\restylefloat{figure}
\usepackage{url}
\usepackage[backref=page]{hyperref}
\hypersetup{colorlinks=true,breaklinks=true,
            linkcolor=darkblue,urlcolor=darkblue,
            anchorcolor=darkblue,citecolor=darkblue}

\usepackage{fancyhdr}
\allowdisplaybreaks[4]

\renewcommand*{\backref}[1]{}
\renewcommand*{\backrefalt}[4]{%
    \ifcase #1 (Not cited.)%
    \or        (Cited on page~#2.)%
    \else      (Cited on pages~#2.)%
    \fi}

\title{An example related to the Erd\H{o}s--Falconer question over arbitrary finite fields}
\author{Brendan Murphy and Giorgis Petridis}
\date{}

\theoremstyle{plain}
\newtheorem{theorem}{Theorem}

\theoremstyle{definition}

\renewcommand*{\backref}[1]{}
\renewcommand*{\backrefalt}[4]{%
    \ifcase #1 (Not cited.)%
    \or        (Cited on page~#2.)%
    \else      (Cited on pages~#2.)%
    \fi}

\newcommand{\F}{\mathbb{F}} 

\begin{document}

\onehalfspacing

\pagenumbering{arabic}

\setcounter{section}{0}

\bibliographystyle{plain}

\maketitle

\begin{abstract}
There exists an infinite family of examples of subsets of $\F_q^2$ with $q^{4/3}$ elements whose distance sets are not the whole of $\F_q$.
\end{abstract}

\let\thefootnote\relax\footnotetext{The first author is supported by the Leverhulme grant RPG 2017-371. The second author is supported by the NSF Award 1723016 and gratefully acknowledges the support from the RTG in Algebraic Geometry, Algebra, and Number Theory at the University of Georgia, and from the NSF RTG grant DMS-1344994.}

\section[Introduction]{Introduction}
\label{Introduction}

Iosevich and Rudnev were the first to investigate in \cite{Iosevich-Rudnev2007} the Erd\H{os}--Falconer distance set question (\cite{Erdos1946, Falconer1985}) over finite fields. In this context distance is defined algebraically and the \emph{distance set} of a set $E$ is
\[
\Delta(E) = \{ (x-y) \cdot (x-y) : x,y \in E\}.
\]
The \emph{distance number} of $E$ is $|\Delta(E)|$. Iosevich and Rudnev studied large subsets of vector spaces over finite fields and sought to find a condition on the cardinality of the subset that guarantees that its distance set is the entire finite field. They proved the following lower bound ($\F_q$ denotes a finite field of $q$ elements).

\begin{theorem}[Iosevich--Rudnev] \label{IR}
For all $d \geq 2$ and all odd $q$, if $E \subseteq \F_q^d$ satisfies $|E| > 4 q^{(d+1)/2}$, then $\Delta(E) = \F_q$.  
\end{theorem}

The exponent $(d+1)/2$ cannot be improved for odd dimensions $d$. Hart, Iosevich, Koh, and Rudnev proved in \cite{HIKR2011} that there exists an absolute $c>0$ and subsets $E \subset \F_q^d$ for any odd $d$, such that $|E| \geq c q^{(d+1)/2}$ and $\Delta(E) \neq \F_q$. Theorem~\ref{IR} remains the state-of-the-art for even $d$.

In this short note, we construct an infinite family of examples for $d =2$, which shows that the exponent in Theorem~\ref{IR} must be at least as big as $4/3$. Our examples disprove Conjecture 1.1 in~\cite{Iosevich-Rudnev2007} for $d=2$.

\begin{theorem} \label{main}
Let $p$ be a prime. There exist infinitely many powers $q$ of $p$ and a subset $E \subset \F_q^2$ such that $|E| = q^{4/3}$ and $\Delta(E) \neq \F_q$. In fact, as $q \to \infty$, $|\Delta(E)| / q \leq 1/2$.
\end{theorem}

The construction is very similar to a construction in~\cite{GPPDA}. Some remarks.
\begin{enumerate}
\item The examples in Theorem~\ref{main} stem from properties of vector spaces over subfields. This means that they cannot be extended to prime order fields. They show that vector spaces over subfields (and not only subfields) can obstruct growth in geometric questions.
\item The examples in Theorem~\ref{main} complement results in~\cite{CEHIK2012,HLRN2016}, which show that if $E \subseteq \F_q^2$ and $|E| > q^{4/3}$, then both the distance number and, the intuitively defined, pinned distance number of $E$ is a constant multiple of $q$. Theorem~\ref{main} shows that for such $E$ $\Delta(E)$ does not always contain, say, $99\%$ of the elements of $\F_q$.
\item Theorem~\ref{main} does not rule out the existence of a positive $\delta>0$ with the property that, for all sufficiently large $q$, any $E \subseteq \F_q^2$ with $|E| > q^{4/3 - \delta}$ elements satisfies $|\Delta(E)| > q/2$. A result along these lines was proved in~\cite{GPPDA} for an analogous question. 
\end{enumerate}

\section[Proof of Theorem~2]{Proof of Theorem~\ref{main}}

We prove Theorem~\ref{main} for odd primes $p$ that are congruent to 3 modulo 4. The cases where $p=2$ or $p$ is congruent to $1$ modulo 4 are similar.

We start with $\F_p$ and adjoin a square root of -1, denoted by  $i$, to get the quadratic extension $\F_p[i]$. Let $F$ be any finite extension of $\F_p[i]$ and $\F_q$ be any cubic extension of $F$. In other words we may take $q = p^{6r}$ for every positive integer $r$. 

Consider $\F_q$ as a vector space over $F$ and let $V$ be any 2-dimensional vector subspace. It is shown in \cite[Proposition~11]{GPPDA} that $VV \neq \F_q$ and that $|VV|/q \to 1/2$ as ($r$ and hence) $q\to \infty$.

We are now in position to describe $E$. Set
\[
E = \{ (u, iv) : u,v \in V\} \subset \F_q^2.
\]
Finally note that, because of closure of $V$ under addition and subtraction,
\begin{align*}
\Delta(E) 
& = \{ (u_1-u_2)^2 + (iv_1 - iv_2)^2 : u_1,u_2, v_1,v_2 \in V \} \\
& = \{ u^2  - v^2 : u,v \in V \} \\
& = \{ (u-v)(u+v) : u,v \in V \} \\
& \subseteq 
VV.
\end{align*}
Therefore $\Delta(E) \neq \F_q$ while $|E| = |V|^2 = |F|^4 = q^{4/3}$. \qedhere 

Note that when $q$ is odd we have $\Delta(E) = VV$.

\phantomsection

\addcontentsline{toc}{section}{References}

\bibliography{all}

\bigskip 

\hspace{20pt} Department of Mathematics, University of Bristol, and Heilbronn Institute of Mathematical Research, Bristol, UK.

\hspace{20pt} Department of Mathematics, University of Georgia, Athens, GA 30602, USA.

\hspace{20pt} \textit{Email addresses}: \href{mailto:brendan.murphy@bristol.ac.uk}{brendan.murphy@bristol.ac.uk} and  \href{mailto:giorgis@cantab.net}{giorgis@cantab.net}.

\end{document}